\documentclass[12pt]{amsart}
\usepackage{color}
\usepackage{epsfig}
\usepackage{verbatim}

\begin{document}
\title {Applications of the Strong Splitter Theorem: decomposition results}

\maketitle 
\begin {center}
S. R. Kingan 
\footnote{The author is partially supported by  PSC-CUNY grant number 66305-00 44.} \\     
Department of Mathematics \\
Brooklyn College, City University of New York\\
 Brooklyn, NY 11210\\
skingan@brooklyn.cuny.edu\\  
\end {center}
\bigskip

\begin{abstract} We use the Strong Splitter Theorem to decompose the excluded minor class of binary matroids with no $E_4$-minor. Using this theorem we can get the 3-decomposers and the extremal internally 4-connected matroids as well as any other important matroids in the class. The matroid $E_4$ is a self-dual 10-element binary 3-connected matroid that plays a useful role in structural results. It is a single-element coextension of $P_9$, which is a single-element extension of the 4-wheel. We show that the extremal matroids in this class are the binary rank-$r$ spikes $Z_r$, the rank 3 and 4 projective geometries $F_7$ and $PG(3,2)$, respectively, the 17-element internally 4-connected matroid $R_{17}$, and one 12-element rank-6 matroid. All the other 3-connected members have $P_9$ or $P_9^*$ as 3-decomposers. As immediate corollaries we get decomposition results for $EX[P_9^*]$ and $EX[P_9]$ as well as the internally 4-connected members of these classes. 
\end{abstract}

\bigskip 

\section {\bf Introduction}

 The Splitter Theorem states that, if $N$ is a $3$-connected proper minor of a 3-connected matroid $M$ such that, if $N$ is a wheel or whirl then $M$ has no larger wheel or whirl, respectively, then there is a sequence $M_0, \dots , M_n$ of $3$-connected matroids with $M_0\cong N$, $M_n=M$, and for $i\in \{1, \dots , n\}$, $M_i$ is a single-element extension or coextension of $M_{i-1}$ [\ref{Seymour1980}]. The Strong Splitter Theorem imposes an ordering on the way $N$ can be extended and coextended. We can obtain up to isomorphism $M$ starting with $N$ and at each step doing a 3-connected single-element extension or coextension, such that at most two consecutive single-element extensions occur in the sequence (unless the rank of the matroids involved is $r(M)$). Moreover, if two consecutive single-element extensions by elements $\{e_1, e_2\}$ are followed by a coextension by element $f$, then $\{e_1, e_2, f\}$ form a triad in the resulting matroid [\ref{KinganLemos2014}]. The Strong Splitter Theorem gives the most efficient way of generating a class of 3-connected matroids, thereby reducing computations.

Let $EX[M_1\dots, M_k]$ denote the class of 3-connected binary matroids with no minors isomorphic to $M_1, \dots, M_k$.  In this paper we demonstrate the usefulness of the Strong Splitter Theorem by characterizing an excluded minor class. The general strategy in our characterization is to find 3-decomposers and the extremal 3-connected matroids that cannot be decomposed. Note that it is sufficient to focus on the 3-connected members of an excluded minor class, since matroids that are not 3-connected can be pieced together from 3-connected matroids using the operations of 1-sum and 2-sum. Moreover, when we exclude a non-regular matroids all the regular matroids are automatically in the class. Thus, in the statement of our theorems we identify only the 3-connected non-regular matroids.

A matrix representation for $P_9$ is given below. It is a single-element extension of the 4-wheel. The matroid $P_9^*$ is its dual and the matroid $E_4$ is a single-element extension of $P_9^*$.

\[ 
P_{9}=\left[ 
\begin{array}{c|ccccc}
&   0&1&1&1&1 \\
I_4&1&0&1&1&1 \\
&   1&1&0&1&0 \\
&   1&1&1&1&0
\end{array} 
\right] 
P_9^*=\left[ 
\begin{array}{c|cccc}
&    0&1&1&1 \\
&    1&0&1&1 \\
I_5& 1&1&0&1 \\
&    1&1&1&1 \\
&    1&1&0&0
\end{array} 
\right] 
E_4=\left[ 
\begin{array}{c|ccccc}
&    0&1&1&1&1 \\
&    1&0&1&1&0 \\
I_5& 1&1&0&1&0 \\
&    1&1&1&1&0 \\
&    1&1&0&0&1
\end{array} 
\right] 
\] 

The main theorem in this paper characterizes $EX[E_4]$. As immediate corrollaries we get characterizations for $EX[P_9^*]$ and $EX[P_9]$. 
Observe that $F_7$  and $PG(3,2)$ shown below are in $EX[E_4]$ as well as in $EX[P_9^*]$ because they have rank 3 and 4, respectively, and $E_4$ and $P_9^*$ have rank 5. 

In [\ref{Oxley1987}] Oxley proved that a $3$-connected binary non-regular matroid $M$ has no minor isomorphic to $P_9$ or $P_9^*$ if and only if $M$ is isomorphic to $F_7$, $F_7^*$, $Z_r$, $Z_r^*$, $Z_r\backslash b_r$, or  $Z_r\backslash c_r$,  for some $r\ge 4$. The matroid $Z_r$ is represented by the matrix $[I_r | D]$ where $D$ has $r+1$ columns labeled $b_1 \dots b_r, c_r$ with zeros on the diagonal and ones elsewhere. The matroids $Z_r\backslash b_r$ and $Z_r\backslash c_r$ are its 3-connected single-element deletion-minors and $Z_r\backslash \{b_r, c_r\}\cong Z_{r-1}^*$. Moreover, $Z_r\backslash b_r$ and $Z_r\backslash c_r$  are self-dual. Thus, the infinite family $Z_r$ for $r\ge 4$ and its 3-connected deletion minors are in $EX[E_4]$ and $EX[P_9^*]$.

\[
F_7=\left[ 
\begin{array}{c|cccc}
&   0&1&1&1   \\
I_3&1&0&1&1  \\
&   1&1&0&1  \\
\end{array} 
\right] 
PG(3, 2)=\left[ 
\begin{array}{c|cccccccccccc}
&   0&0&0&0&1&1&1&1&1&1&1  \\
I_4& 0&1&1&1&0&0&0&1&1&1&1 \\
&   1&0&1&1&0&1&1&0&0&1&1 \\
&   1&1&0&1&1&0&1&0&1&0&1
\end{array} 
\right] 
\] 

In addition to these matroids, we prove that $EX[E_4]$ has the 17-element rank-5 internally 4-connected matroid $R_{17}$. This matroid appears in [\ref{KinganLemos2014}] as the extremal matroid for $EX[(K_5\backslash e)^*]$. The class also has a 12-element rank-6 3-connected matroid $M_{12}$ that is a splitter (every 3-connected single-element extension and coextension is not in the class). This matroid $M_{12}$ is not internally 4-connected. The rest of the 3-connected matroids have a non-minimal exact 3-separation induced by $P_9$ or $P_9^*$. Matrix representations for $M_{12}$ and $R_{17}$ are shown below.

\[
M_{12}=\left[ 
\begin{array}{c|cccccc}
&     0 &  1  & 1 &  1 &  1 &  1  \\
&     1 &  0  & 1 &  1 &  0 &  0 \\
I_6&  1 &  1  & 0 &  1 &  1 &  0 \\
&     1 &  1  & 1 &  1 &  0 &  1  \\
&     1 &  1  & 0 &  0 &  0 &  1 \\
&     1 &  0  & 0 &  1 &  1 &  1 
\end{array}
\right]  
R_{17}=\left[ 
\begin{array}{c|cccccccccccc}
&    1&0&0&1&1&0&0&1&1&1&1&1 \\
&    1&1&0&0&1&1&1&0&0&1&1&1 \\
I_5& 1&1&1&0&0&0&1&1&1&0&1&1 \\
&    0&1&1&1&0&1&0&0&1&1&1&1 \\
&    0&0&1&1&1&1&1&1&0&0&1&0
\end{array}
\right] 
\]

\noindent The next result is the main theorem in this paper.
\bigskip

\noindent{\bf Theorem 1.1.} {\it Suppose $M$ is a binary $3$-connected  non-regular matroid with no $E_4$-minor. Then either $P_9$ or $P_9^*$ are $3$-decomposers for $M$ or $M$ or $M^*$ is isomorphic to a $3$-connected deletion-minor of $Z_r$ for $r\ge 4$, $F_7$, $PG(3,2)$, $R_{17}$, or $M_{12}$. }
\bigskip 

As an immediate corollary we can characterize $EX[P_9^*]$. The rank-5 internally 4-connected extremal matroid in $EX[P_9^*]$ is $R_{17}\backslash \{17\}$. For convenience of notation we call it $R_{16}$.  
\bigskip

\noindent{\bf Corollary 1.2} {\it Suppose $M$ is a binary $3$-connected  non-regular matroid with no $P_9^*$-minor. Then one of the following holds:
\begin {enumerate}
\item [(i)] $M$ is isomorphic to $F_7$, $F_7^*$, $Z_r$, $Z_r^*$, $Z_r\backslash b_r$, or $Z_r\backslash c_r$, for $r\ge 4$;
\item [(ii)]  $P_9$ is a $3$-decomposer for $M$; or 
\item [(iii)] $M$ is isomorphic to a $3$-connected deletion-minor of $PG(3,2)$ or $R_{16}$.
\end {enumerate}}
\bigskip

\noindent By duality we get the following characterization for $EX[P_9]$. 
\bigskip

\noindent{\bf Corollary 1.3} {\it Suppose $M$ is a binary $3$-connected  non-regular matroid with no $P_9$-minor. Then one of the following holds:
\begin {enumerate}
\item [(i)] $M$ is isomorphic to $F_7$, $F_7^*$, $Z_r$, $Z_r^*$, $Z_r\backslash b_r$, or $Z_r\backslash c_r$, for $r\ge 4$;
\item [(ii)]  $P_9^*$ is a $3$-decomposer for $M$; or 
\item [(iii)] $M$ is isomorphic to a $3$-connected contraction-minor of $R_{16}$.
\end {enumerate}}
\bigskip

We end this section by describing our method for calculating extensions and coextensions. This method along with the Strong Splitter Theorem leads to relatively short proofs. 

Let $N$ be a $GF(q)$-representable $n$-element rank-$r$  matroid  represented by the matrix $A=[I_r|D]$ over $GF(q)$. The columns of A may be viewed as a subset of the columns of the matrix that represents the projective geometry $PG(r - 1, q)$. Let $M$ be a simple single-element extension of $N$ over $GF(q)$. Then $N=M\backslash e$ and $M$ may be represented by $[I_r|D']$, where $D'$ is the same as $D$, but with one additional column corresponding to the element $e$. The new column is distinct from the existing columns and has at least two non-zero elements. If the existing columns are labeled $\{1, \dots , r, \dots , n\}$, then the new column is labeled $(n+1)$.

Suppose $M$ is a cosimple single-element coextension of $N$ over $GF(q)$. Then $N=M/f$ and $M$ may be represented by the matrix $[I_{r+1}| D'']$, where $D''$ is the same as $D$, but with one additional row. The new row is distinct from the existing rows and has at least two non-zero elements. The columns of $[I_{r+1}| D'']$ are labeled $\{1, \dots , r+1, r+2, \dots , n, n+1\}$. The coextension element $f$ corresponds to column $r+1$. The coextension row is selected from $PG(n-r, q)$. We can visualize the new element $f$  as appearing in the new dimension and lifting several points into the higher dimension. Observe that $f$ forms a cocircuit with the elements corresponding to the non-zero elements in the new row.  Note that in $[I_{r+1}|D'']$ the labels of columns beyond $r$ are increased by 1 to accomodate the new column $r+1$. 

We refer to the simple single-element extensions of $N$ as Type (i) matroids and the cosimple single-element coextensions of $N$ as Type (ii) matroids. The structure of Type (i) and Type (ii) matroids are shown in Figure 1.

\begin{figure}[h]
\centering
\epsfxsize 5in \epsfbox{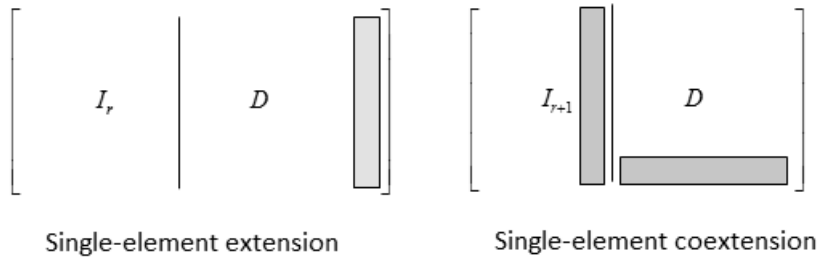}
\caption{Structure of Type (i) and Type (ii) matroids }
\end{figure}

Once the simple single-element extensions (Type (i) matroids) and cosimple single-element coextensions (Type (ii) matroids) are determined, the number of permissable rows and columns give a bound on the choices for the cosimple single-element extensions of the Type (i) matroids and the simple single-element extensions of the Type (ii) matroids, respectively. 

The structure of the cosimple single-element coextensions of a Type (i) matroid and the simple single-element extensions of a Type (ii) matroid are shown in Figure 2.
 
\begin{figure}[h]
\centering
\epsfxsize 5in \epsfbox{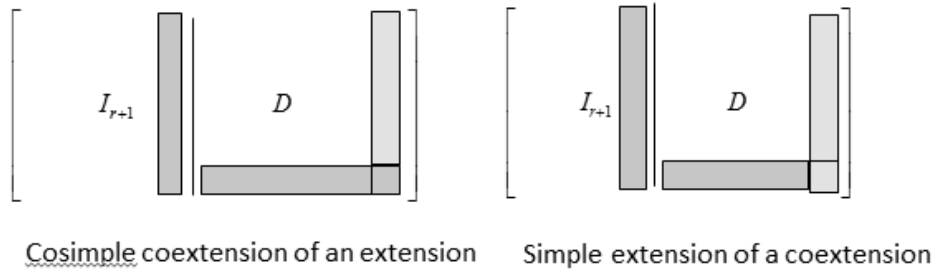}
\caption{Structure of $M$, where $|E(M)-E(N)|=2$ }
\end{figure}
 
\noindent When computing the cosimple single-element coextension of a Type (i) matroid, there are three types of rows that may be inserted into the last row.

\begin {enumerate}
\item[(I)] rows that can be added to $N$ to obtain a coextension with a 0 or 1 as the last entry (or as many as the entries in $GF(q)$ for higher order fields); 
\item[(II)] the identity rows with a 1 in the last position; and
 \item[(III)] rows ``in-series" to the right-hand side of the matrix with the last entry reversed. 
\end{enumerate}

\noindent When computing the simple single-element extension of a Type (ii) matroid, there are three types of rows that may be inserted into the last column.

\begin {enumerate}
\item[(I)] columns that can be added to $N$ to obtain an extension with a 0 or 1 as the last entry (or as many as the entries in $GF(q)$ for higher order fields); 
\item[(II)] the identity columns with a 1 in the last position; and
 \item[(III)] columns ``in-parallel" to the right-hand side of matrix with the last entry reversed. 
\end{enumerate}

Suppose $N'$ is a simple double-element extension of $N$ formed by adding columns $e_1$ and $e_2$ and $M$ is a cosimple single-element coextension of $N'$ by element $f$. Then, by the Strong Splitter Theorem  $M\backslash e_1$ or $M\backslash e_2$ is 3-connected except when $\{e_1, e_2, f\}$ is a triad. Thus, the only coextension of $N'$ we must check is the one formed by adding row $[0 0 \dots 0 1 1]$. Additionally, there is no need to calculate the cosimple single-element extensions of the $k$-element extensions of $N$, for $k\ge 3$. This greatly reduces computations. 

\section {Characterization of $EX[E_4]$}

All proofs of excluded minor characterizations of binary non-regular matroids begin the same way. Tutte proved that a binary matroid is non-regular if and only if it has no minor isomorphic to $F_7$ or $F_7^*$ [\ref{Oxley2012}, 10.1.2].  Observe that $F_7=PG(2,2)$ and as such has no extensions in the class of binary matroids. Coextensions of $F_7$ are duals of extensions of $F_7^*$. Thus we may focus on the extensions of $F_7^*$.  Observe that $AG(3, 2)$ and $S_{8}$ are the two non-isomorphic 3-connected single-element extensions of $F_7^*$. Since they are self-dual, they are also the coextensions of $F_7$. The matroid $S_8$ has two non-isomorphic 3-connected single-element extensions $P_9$ and $Z_4$ and $AG(3, 2)$ has one 3-connected single-element extension $Z_4$. The simple extensions and cosimple-single-element coextensions of $P_9$ are given in Table 1a and 1b.

\small
 \begin{center}
\begin{tabular}{|c|p{15em}|c|}
\hline
\bf{Matroid}& \bf{Extension Columns} & {\bf Name}   \\ \hline  
$P_9$ &  $\bf [1 1 1 0]$  & $D_1$   \\ \hline
&  $\bf [1 0 0 1]$ $[0 1 0 1]$ $[0 1 1 0]$, $[1 0 1 0]$  & $D_2$               \\ \hline
& $\bf [0 0 1 1]$    &   $D_3$          \\  \hline \hline
$D_1$ & $[0 1 0 1]$ $[0 1 1 0]$ $[1 0 0 1]$ $[1 0 1 0]$ & $X_1$ \\  \hline  
& $\bf [0 0 1 1]$ & $X_2$   \\  \hline 
$D_2$ & $\bf [1 0 1 0]$ $[1 1 1 0]$  & $X_1$   \\  \hline
& $\bf [0 0 1 1 ]$ $[0 1 0 1]$ $[0 1 1 0]$ & $X_3$   \\  \hline 
$D_3$ & $[1 1 1 0]$  & $X_2$   \\  \hline
& $ [0 1 0 1]$ $[0 1 1 0]$ $[1 0 0 1]$ $[1 0 1 0]$ & $X_3$ \\ \hline 
$X_1$ & $\bf [0 0 1 1]$ $[0 1 0 1]$ $[0 1 1 0]$  & $Y_1$   \\  \hline
& $\bf [1 1 1 0]$ & $Y_2$ \\ \hline
$X_2$ & $[0 1 0 1]$ $[0 1 1 0]$ $[1 0 0 1]$ $[1 0 1 0]$  & $Y_1$ \\ \hline 
$X_3$ & $ [0 1 0 1]$ $[0 1 1 0]$ $[1 0 1 0]$ $[1 1 1 0]$  & $Y_1$   \\  \hline 
\end{tabular}
 \end{center}
 
 \begin{center} Table 1a: Rank 4 extensions of $P_9$ \end{center}

 \begin{center}
\begin{tabular}{|p{28em}|c|}
\hline
\bf{Coextension Rows} & {\bf Name}     \\   \hline \hline
$[1 1 0 0 0]$ $[1 1 1 1 1]$    & $E_1$  \\  \hline
$[1 1 0 1 1]$ $[1 1 1 0 0]$    & $E_2$  \\  \hline
$[1 1 0 0 1]$ $[1 1 1 0 1]$    &  $E_3$  \\  \hline
$[0 1 0 0 1]$ $[0 1 0 1 0]$ $[0 1 1 0 1]$ $[0 1 1 1 0]$ $[1 0 0 0 1]$ $[1 0 0 1 0]$ $[1 0 1 0 1]$ $[1 0 1 1 0]$    &  $E_4$  \\  \hline
$[0 1 0 1 1]$ $[0 1 1 0 0]$ $[1 0 0 1 1]$ $[1 0 1 0 0]$   &  $E_5$  \\  \hline
$[0 0 1 0 1]$ $[0 0 1 1 0]$    &  $E_6$  \\  \hline
$[0 0 1 1 1]$                  &  $E_6^*$  \\  \hline
$[0 0 0 1 1]$      &  $E_7$  \\  \hline
\end{tabular}
 \end{center}
 
 \begin{center} Table 1b: Single-element coextensions of $P_9$ \end{center}
\normalsize

Suppose $M$ is a $3$-connected binary non-regular matroid with a $P_9$-minor. From Tables 1a and 1b we see that $P_9$ has three non-isomorphic simple single-element extensions, $D_1$, $D_2$, and $D_3$, and eight non-isomorphic cosimple single-element coextensions $E_1$, $E_2$, $E_3$, $E_4$, $E_5$, $E_6$, $E_6^*$, and $E_7$.

\tiny
\[ 
D_1=\left[ 
\begin{array}{c|cccccc}
&   0&1&1&1&1&1 \\
I_4&1&0&1&1&1&1 \\
&   1&1&0&1&0&1 \\
&   1&1&1&1&0&0
\end{array} 
\right] 
D_2=\left[ 
\begin{array}{c|cccccc}
&   0&1&1&1&1&1 \\
I_4&1&0&1&1&1&0 \\
&   1&1&0&1&0&0 \\
&   1&1&1&1&0&1
\end{array} 
\right] 
D_3=\left[ 
\begin{array}{c|cccccc}
&   0&1&1&1&1&0 \\
I_4&1&0&1&1&1&0 \\
&   1&1&0&1&0&1 \\
&   1&1&1&1&0&1
\end{array} 
\right] 
\] 
\normalsize

\noindent 
\tiny
\[ 
E_1=\left[ 
\begin{array}{c|ccccc}
&    0&1&1&1&1 \\
&    1&0&1&1&1 \\
I_5& 1&1&0&1&0 \\
&    1&1&1&1&0 \\
&    1&1&0&0&0
\end{array} 
\right] 
E_2=\left[ 
\begin{array}{c|ccccc}
&    0&1&1&1&1 \\
&    1&0&1&1&1 \\
I_5& 1&1&0&1&0 \\
&    1&1&1&1&0 \\
&    1&1&0&1&1
\end{array} 
\right] 
E_3=\left[ 
\begin{array}{c|ccccc}
&    0&1&1&1&1 \\
&    1&0&1&1&1 \\
I_5& 1&1&0&1&0 \\
&    1&1&1&1&0 \\
&    1&1&0&0&1
\end{array} 
\right] 
E_4=\left[ 
\begin{array}{c|ccccc}
&    0&1&1&1&1 \\
&    1&0&1&1&1 \\
I_5& 1&1&0&1&0 \\
&    1&1&1&1&0 \\
&    0&1&0&0&1
\end{array} 
\right] 
\] 
\normalsize

\tiny
\[ 
E_5=\left[ 
\begin{array}{c|ccccc}
&    0&1&1&1&1 \\
&    1&0&1&1&1 \\
I_5& 1&1&0&1&0 \\
&    1&1&1&1&0 \\
&    1&0&1&0&0
\end{array} 
\right] 
E_6=\left[ 
\begin{array}{c|ccccc}
&    0&1&1&1&1 \\
&    1&0&1&1&1 \\
I_5& 1&1&0&1&0 \\
&    1&1&1&1&0 \\
&    0&0&1&0&1
\end{array} 
\right] 
E_6^*=\left[ 
\begin{array}{c|ccccc}
&    0&1&1&1&1 \\
&    1&0&1&1&1 \\
I_5& 1&1&0&1&0 \\
&    1&1&1&1&0 \\
&    0&0&1&1&1
\end{array} 
\right] 
E_7=\left[ 
\begin{array}{c|ccccc}
&    0&1&1&1&1 \\
&    1&0&1&1&1 \\
I_5& 1&1&0&1&0 \\
&    1&1&1&1&0 \\
&    0&0&0&1&1
\end{array} 
\right] 
\] 
\normalsize 

\bigskip
\noindent {\bf Lemma 2.1.} {\it  If $M$ has a $P_9$- or $P_9^*$-minor, but no $D_2$, $D_2^*$, $E_4$, or $E_5$-minor, then $P_9$ or $P_9^*$ are $3$-decomposers for $M$.}
\bigskip

\noindent {\bf Proof.} The proof of this lemma appears in [\ref{KinganLemos2014}, Theorem 3.1] and is repeated here for convenience. Observe that $P_9$ has  a non-minimal exact $3$-separation $(A, B)$, where $A=\{1, 2, 5, 6\}$ is both a circuit and a cocircuit.   It is easy to check that the set $A=\{1, 2, 5, 6\}$ is both a circuit and a cocircuit in $D_1$ and $D_3$ (note that every column is checked) whereas $D_2$ is internally 4-connected. 
The set $A=\{1, 2, 5, 6\}$ corresponds to $A'=\{1, 2, 6, 7\}$ in the coextension since the fifth column is the coextended element. It can be checked that $\{1, 2, 6, 7\}$ is both a circuit and a cocircuit in $E_1$, $E_2$, $E_3$, $E_6$, $E_6^*$, and $E_7$ (every row is checked disregarding isomorphism).  Further note that  $E_4$ and $E_5$ are self-dual. Thus if $M$ has a $P_9$-minor, but no $D_2$, $D_2^*$, $E_4$, or $E_5$-minor, then $P_9$ or $P_9^*$ is a $3$-decomposer for $M$ using Mayhew, Royle, and Whittle's sufficient one-element check [\ref{MayhewRoyleWhittle2011}]. $\qed$
\bigskip

\noindent {\bf Proof of Theorem 1.1.} Suppose $M$ is a 3-connected binary matroid with no $E_4$-minor. Lemma 2.1 implies that if $M$ has a $P_9$- or $P_9^*$-minor, but no $D_2$, $D_2^*$, $E_4$, or $E_5$-minor, then $P_9$ or $P_9^*$ are $3$-decomposers for $M$. So we must consider matroids that have an $E_5$, $D_2$, or  $D_2^*$-minor, but no $E_4$-minor. 
Consider the 3-connected single-element extensions and coextensions of $E_5$ shown in Tables 2a and 2b.  

\small
 \begin{center}
\begin{tabular}{|p{15em}|c|c|}
\hline
 \bf{Extension Columns} & {\bf Name} & {\bf $E_4$-minor}  \\  \hline \hline
  $[0 0 1 0 1]$ $[0 0 1 1 0]$ $[0 1 0 1 1]$ $[0 1 1 0 0]$  & $A$ &  No   \\  \hline

 $[1 0 0 1 1]$ & $B$   & No   \\  \hline

 $[1 1 0 0 1]$ $[1 1 1 0 1]$   &   $C$   & No   \\  \hline

 $[0 0 0 1 1 ]$ $[0 0 1 1 1]$ $[0 1 0 0 1]$ $[0 1 1 0 1]$  && Yes   \\  \hline

 $0 1 0 1 0]$ $[0 1 1 1 0]$ && Yes     \\  \hline

 $[1 0 0 0 1]$ $[1 0 0 1 0]$  $[1 1 0 1 1]$ $[1 1 1 0 0]$  && Yes   \\  \hline

 $[1 0 1 0 1]$ $[1 0 1 1 0]$ $[1 1 0 0 0]$ $[1 1 1 1 1]$   & & Yes   \\  \hline
\end{tabular}
 \end{center}
 
 \begin{center}  Table 2a: Simple single-element extensions of $E_5$  \end{center}

 \begin{center}
\begin{tabular}{|p{15em}|c|c|c|}
\hline
\bf{Coextension Rows} & {\bf Name} \\      \hline \hline
$[0 0 1 1 1]$ $[0 1 0 0 1]$ $[0 1 0 1 0]$ $[0 1 1 0 0]$  & $A^*$  \\  \hline

$[1 0 0 1 1]$ & $B^*$      \\  \hline

$[1 0 1 0 1]$ $[1 1 1 0 1]$   &   $C^*$    \\  \hline

$[0 0 0 1 1 ]$ $[0 0 1 0 1]$ $[0 1 0 1 1]$ $[0 1 1 0 1]$   & \\  \hline

$[0 0 1 1 0]$ $[0 1 1 1 0]$   &   \\  \hline

$[1 0 0 0 1]$ $[1 0 0 1 0]$  $[1 0 1 1 1]$ $[1 1 1 0 0]$ & \\  \hline

$[1 0 1 0 0]$ $[1 1 0 0 1]$ $[1 1 0 1 0]$ $[1 1 1 1 1]$ &   \\  \hline
\end{tabular}
 \end{center}

 \begin{center}   Table 2b: Cosimple single-element coextensions of $E_5$ \end{center} 
\normalsize

\noindent Observe that $E_5$ has seven non-isomorphic simple single-element extensions all of which have an $E_4$-minor except $A$, $B$ and $C$. Matrix representations for $A$, $B$, and $C$ are given below.

\small
\[ 
A=\left[ 
\begin{array}{c|cccccc}
&    0&1&1&1&1&0 \\
&    1&0&1&1&0&0 \\
I_5& 1&1&0&1&1&1 \\
&    1&1&1&1&0&0 \\
&    1&1&0&0&0&1
\end{array} 
\right] 
B=\left[ 
\begin{array}{c|ccccccc}
&    0&1&1&1&1&1 \\
&    1&0&1&1&0&0 \\
I_5& 1&1&0&1&1&0 \\
&    1&1&1&1&0&1 \\
&    1&1&0&0&0&1
\end{array} 
\right] 
C=\left[ 
\begin{array}{c|ccccccc}
&    0&1&1&1&1&1 \\
&    1&0&1&1&0&1 \\
I_5& 1&1&0&1&1&0 \\
&    1&1&1&1&0&0 \\
&    1&1&0&0&0&1
\end{array} 
\right] 
\] 
\normalsize

\noindent The next lemma is a key lemma for Theorem 1.1. It is also a useful stand-alone result for future structure theorems since $E_4$ and $E_5$ are important binary matroids. For this lemma we use the following representation of $E_5$:

\small
\[ 
E_5=\left[ 
\begin{array}{c|ccccc}
&    0&1&1&1&1 \\
&    1&0&1&1&0 \\
I_5& 1&1&0&1&1 \\
&    1&1&1&1&0 \\
&    1&1&0&0&0
\end{array} 
\right]
\] 
\normalsize

\bigskip

\noindent {\bf Lemma 2.2.} {\it Suppose $M$ is a binary $3$-connected matroid with an $E_5$-minor and no $E_4$-minor.  Then either $M\cong M_{12}$ or $M$ or $M^*$ is isomorphic to $R_{17}$ or its $3$-connected deletion-minors having an $E_5$-minor.}  
\bigskip

\noindent {\bf Proof.} The proof is in three stages. First, we will show that all the coextensions of $A$, $B$, and $C$ have an  $E_4$-minor with the exception of $M_{12}$. Suppose $M$ is a coextension of $A$, $B$, $C$. The three types of rows that may be added to $A$, $B$ and $C$ to obtain $M$ are:

\begin {enumerate}
\item[(I)] rows that can be added to $E_5$ to obtain a coextension with no $E_4$-minor, with a 0 or 1 as the last entry; 
\item[(II)] the identity rows with a 1 in the last position; and
 \item[(III)] the rows ``in-series" to the right-hand side of matrices $A$, $B$, $C$ with the last entry reversed. 
\end{enumerate}

\noindent Type I rows are $[0 0 1 1 1 0]$, $[0 0 1 1 1 1]$ $[0 1 0 0 1 0]$, $[0 1 0 0 1 1]$, $[0 1 0 1 0 0]$, $[0 1 0 1 0 1]$, $[0 1 1 0 0 0]$, $[0 1 1 0 0 1]$, $[1 0 0 1 1 0]$, $[1 0 0 1 1 1]$, $[1 0 1 0 1 0]$, $[1 0 1 0 1 1]$, $[1 1 1 0 1 0]$, and $[1 1 1 0 1 1]$. They are obtained from Table 2b.
Type II rows are $[1 0 0 0 0 1]$, $[0 1 0 0 0 1]$, $[0 0 1 0 0 1]$, $[0 0 0 1 0 1]$, and $[0 0 0 0 1 1]$. 
Type III rows are specific to the matrices $A$, $B$, $C$. For matrix $A$ they are $[0 1 1 1 1 1]$, $[1 0 1 1 0 1]$, $[1 1 0 1 1 0]$, $[1 1 1 1 0 1]$, $[1 1 0 0 0 0]$. For matrix $B$ they are $[0 1 1 1 1 0]$, $[1 0 1 1 0 1]$, $[1 1 0 1 1 1]$, $[1 1 1 1 0 0]$, and $[1 1 0 0 0 0]$. For $C$ they are $[0 1 1 1 1 0]$, $[1 0 1 1 0 0]$, $[1 1 0 1 1 1]$, $[1 1 1 1 0 1]$, and $1 1 0 0 0 0]$.

Most of the above rows result in matroids that have an $E_4$-minor (see bold-face rows in Table 3). Only a few coextensions must be specifically checked for an $E_4$-minor. They are $(A, coextn 11)$, $(B, coextn 8)$, $(C, coextn 8)$, $(C, coextn 9)$, $(C, coextn 10)$, $(C, coextn 12)$, and $(C, coextn 14)$. 

Observe that 
$(A, coextn 11)/11\backslash 3 \cong E_4$, 
$(C, coextn 8)/12\backslash 2 \cong E_4$, 
$(C, coextn 9)/12\backslash 1 \cong E_4$, 
$(C, coextn 10)/12\backslash 10 \cong E_4$, and 
$(C, coextn 14)/12\backslash 6 \cong E_4$. 
Further, $(B, coextn 8)\cong (C, coextn 12)$ and this matroid does not have an $E_4$-minor. This is the matroid $M_{12}$. 

Second, we must establish that $M_{12}$ is a splitter for $EX[E_4]$. By the Splitter Theorem and the fact that $M_{12}$ is self-dual, we only need to check the single-element coextensions of $M_{12}$. From Table 3 observe that $M_{12}$ as a coextension of $C$ may be obtained by adding exactly one row. Thus, there are no further rows that may be added to form coextensions without an $E_4$-minor. It follows that $M_{12}$ is a splitter for the class of binary matroids with no $E_4$-minor.

Third, we must show that  either $M\cong M_{12}$ or $r(M)\le 5$. To show this we compute the simple single-element extensions of $A$, $B$, and $C$ with no $E_4$-minor. From Table 2a the only columns that can be added to $E_5$ to obtain a matroid with no $E_4$-minor are  $[0 0 1 0 1]$, $[0 0 1 1 0]$, $[0 1 0 1 1]$, $[0 1 1 0 0]$ $[1 0 0 1 1]$, $[1 1 0 0 1]$, $[1 1 1 0 1]$. They give the matroids $D$, $E$, $F$, and $G$ shown below.

\small
\[ 
D=\left[ 
\begin{array}{c|ccccccc}
&    0&1&1&1&1&0&0 \\
&    1&0&1&1&0&0&0 \\
I_5& 1&1&0&1&1&1&1 \\
&    1&1&1&1&0&0&1 \\
&    1&1&0&0&0&1&0
\end{array} 
\right] 
E=\left[ 
\begin{array}{c|ccccccc}
&    0&1&1&1&1&0&0 \\
&    1&0&1&1&0&0&1 \\
I_5& 1&1&0&1&1&1&0 \\
&    1&1&1&1&0&0&1 \\
&    1&1&0&0&0&1&1
\end{array} 
\right] 
\] 

\[ 
F=\left[ 
\begin{array}{c|ccccccc}
&    0&1&1&1&1&0&1 \\
&    1&0&1&1&0&0&1 \\
I_5& 1&1&0&1&1&1&0 \\
&    1&1&1&1&0&0&0 \\
&    1&1&0&0&0&1&1
\end{array} 
\right] 
G=\left[ 
\begin{array}{c|ccccccc}
&    0&1&1&1&1&0&1 \\
&    1&0&1&1&0&0&1 \\
I_5& 1&1&0&1&1&1&1 \\
&    1&1&1&1&0&0&0 \\
&    1&1&0&0&0&1&1
\end{array} 
\right] 
\] 
\normalsize

Specifically, adding to $A$ column $[0 0 1 1 0]$, $[0 1 1 0 0]$, or $[1 0 0 1 1]$ gives $D$; adding column $[0 1 0 1 1]$ gives $E$; adding $[1 1 0 0 1]$ gives $F$; and adding $[1 1 1 0 1]$ gives $G$. Similarly, we can check that $B$ extends to $D$ and $F$ and $C$ extends to $F$ and $G$. Observe that adding all seven columns to $E_5$ gives the 17-element matroid shown below which is isomorphic to the representation of $R_{17}$ shown in the introduction. (Note that our strategy of throwing away all the matroids that can be decomposed and directly obtaining the extremal internally 4-connected matroid side-steps having to calculate all the internally 4-connected matroids.)

\small 
\[
R_{17}=\left[ 
\begin{array}{c|cccccccccccc}
&    0&1&1&1&1&0&0&0&0&1&1&1 \\
&    1&0&1&1&0&0&0&1&1&0&1&1 \\
I_5& 1&1&0&1&1&1&1&0&1&0&0&1 \\
&    1&1&1&1&0&0&1&1&0&1&0&0 \\
&    1&1&0&0&0&1&0&1&0&1&1&1
\end{array}
\right] 
\] 
\normalsize

\tiny
 \begin{center}
\begin{tabular}{|c|c|p{35em}|}
\hline
\bf{Matroid} & Name  & {\bf Coextension Row}  \\  \hline \hline
$A$ & coext 1 & $\bf [0 0 0 0 1 1]$  $\bf [0 0 0 1 0 1]$   $[0 0 1 0 1 0]$ $[0 1 1 0 1 0]$ $[1 0 1 1 1 1]$ $[1 1 1 0 0 1]$               \\  \hline

& coext 2 & $[0 0 0 1 1 0]$  $[1 1 0 0 1 1]$  $[1 1 0 1 0 1]$    \\  \hline

& coext 3 & $[0 0 0 1 1 1]$  $\bf [1 0 1 0 1 1]$ $\bf [1 1 1 0 1 1]$ \\  \hline

& coext 4 & $\bf [0 0 1 0 0 1]$ $[0 1 0 1 1 0]$ $\bf [0 1 1 1 1 1]$       \\ \hline

& coext 5 & $[0 0 1 0 1 1]$ $[0 1 1 0 1 1]$  $\bf [1 0 0 1 1 1]$      \\ \hline

& coext 6 & $[0 0 1 1 0 0]$ $[0 1 1 1 0 0]$  $\bf [1 1 0 0 0 0]$      \\ \hline

& coext 7 & $[0 0 1 1 0 1]$  $\bf [0 1 0 0 1 0]$ $\bf [0 1 0 1 0 0]$   $[0 1 1 1 0 1]$  $[1 0 1 1 1 0]$  $[1 1 1 0 0 0]$     \\ \hline

& coext 8 & $\bf [0 0 1 1 1 0]$ $\bf [0 1 1 0 0 0]$  $\bf [1 0 1 1 0 1]$   $[1 1 0 0 1 0]$  $[1 1 0 1 0 0]$   $\bf [1 1 1 1 0 1]$      \\ \hline

& coext 9 & $\bf [0 0 1 1 1 1]$ $\bf [0 1 1 0 0 1]$   $[1 0 0 0 1 1]$  $[1 0 0 1 0 1]$    $\bf [1 0 1 0 1 0]$  $\bf [1 1 1 0 1 0]$     \\ \hline

& coext 10 &   $\bf [0 1 0 0 0 1]$   $\bf [1 0 0 0 1 0]$  $[1 0 0 1 0 0]$     \\ \hline

& coext 11 &  $\bf [0 1 0 0 1 1]$ $\bf [0 1 0 1 0 1]$ $\bf [1 0 0 1 1 0]$    \\ \hline

& coext 12 & $[0 1 0 1 1 1]$ \\ \hline

& coext 13 & $\bf [1 0 0 0 0 1]$   $[1 0 1 0 0 0]$ $[1 1 1 1 1 0]$    \\ \hline

& coext 14 & $[1 0 1 0 0 1]$  $\bf [1 1 0 1 1 0]$  $[1 1 1 1 1 1]$     \\ \hline

 \hline

$B$ & coext 1  & $\bf [0 0 0 0 1 1]$ $\bf [0 0 0 1 0 1]$  $[0 0 0 1 1 0]$  $\bf [0 0 1 0 0 1]$  $[0 0 1 0 1 0]$  $\bf [0 0 1 1 1 1]$ $\bf [0 1 0 0 1 0]$ $\bf [0 1 0 1 0 0]$   $[0 1 0 1 1 1 ]$   $\bf [0 1 1 0 0 0]$  $[0 1 1 0 1 1]$  $\bf [0 1 1 1 1 0]$             \\  \hline

& coext 2 &  $[0 0 0 1 1 1]$ $[0 0 1 0 1 1]$ $[0 1 0 1 1 0]$ $[0 1 1 0 1 0]$            \\  \hline

& coext 3 &  $[0 0 1 1 0 0]$ $\bf [0 1 0 0 0 1]$    $[0 1 1 1 0 1]$           \\  \hline

& coext 4 &  $[0 0 1 1 0 1]$    $\bf [0 0 1 1 1 0]$ $\bf [0 1 0 0 1 1]$ $\bf [0 1 0 1 0 1]$ $\bf [0 1 1 0 0 1]$   $[0 1 1 1 0 0]$             \\   \hline

& coext 5 &   $\bf [1 0 0 0 0 1]$   $[1 0 0 0 1 0]$ $[1 0 0 1 0 0]$ $[1 0 1 0 0 0]$   $\bf [1 0 1 1 0 1]$   $[1 0 1 1 1 0]$   $\bf [1 1 0 0 0 0]$   $[1 1 0 0 1 1]$ $[1 1 0 1 0 1]$ $[1 1 1 0 0 1]$   $\bf [1 1 1 1 0 0]$   $[1 1 1 1 1 1]$             \\   \hline

& coext 6 & $[1 0 0 0 1 1]$ $[1 0 0 1 0 1]$   $\bf [1 0 1 0 1 0]$    $[1 0 1 1 1 1]$ $[1 1 1 0 0 0]$   $\bf [1 1 1 0 1 1]$              \\   \hline

& coext 7 &  $\bf [1 0 0 1 1 0]$   $[1 0 1 0 0 1]$ $[1 1 0 0 1 0]$ $[1 1 0 1 0 0]$   $\bf [1 1 0 1 1 1]$   $[1 1 1 1 1 0]$            \\   \hline

& coext 8 &   $\bf [1 0 0 1 1 1]$ $\bf [1 0 1 0 1 1]$ $\bf [1 1 1 0 1 0]$             \\   \hline\hline

$C$ & coext 1 & $\bf [0 0 0 0 1 1]$ $\bf [0 0 0 1 0 1]$ $\bf [0 0 1 0 0 1]$  $\bf [0 0 1 1 1 1]$         $\bf [0 1 0 0 1 0]$ $\bf [0 1 0 1 0 0]$  $\bf [0 1 1 0 0 0]$  $\bf [0 1 1 1 1 0]$ \color {black} $[1 0 0 0 1 0]$ $[1 0 0 1 0 0]$ $[1 0 1 0 0 0]$ $[1 0 1 1 1 0]$ $[1 1 0 0 1 1]$ $[1 1 0 1 0 1]$ $[1 1 1 0 0 1]$ $[1 1 1 1 1 1]$                   \\  \hline

& coext 2 & $[0 0 0 1 1 0]$  $[0 1 0 1 1 1]$     \\  \hline

& coext 3 & $[0 0 0 1 1 1]$  $[0 1 0 1 1 0]$   $\bf [1 0 0 1 1 0]$  $\bf [1 1 0 1 1 1]$   \\  \hline

& coext 4 & $[0 0 1 0 1 0]$  $[0 1 1 0 1 1]$    \\  \hline

& coext 5 & $[0 0 1 0 1 1]$  $[0 1 1 0 1 0]$  $\bf [1 0 1 0 1 0]$ $\bf [1 1 1 0 1 1]$     \\  \hline

& coext 6 & $[0 0 1 1 0 0]$  $[0 1 1 1 0 1]$     \\  \hline

& coext 7 & $[0 0 1 1 0 1]$  $[0 1 1 1 0 0]$   $\bf [1 0 1 1 0 0]$ $\bf [1 1 1 1 0 1]$    \\  \hline

& coext 8 &  $\bf [0 0 1 1 1 0]$  $\bf [0 1 0 0 1 1]$ $\bf [0 1 0 1 0 1]$ $\bf [0 1 1 0 0 1]$     \\  \hline

& coext 9 &   $\bf [0 1 0 0 0 1]$     \\  \hline

& coext 10 & $\bf [1 0 0 0 0 1]$  $\bf [1 1 0 0 0 0]$   \\  \hline

& coext 11 & $[1 0 0 0 1 1]$  $[1 0 0 1 0 1]$ $[1 0 1 1 1 1]$ $[1 1 1 0 0 0]$    \\  \hline

& coext 12 &    $\bf [1 0 0 1 1 1]$     \\  \hline

& coext 13 & $[1 0 1 0 0 1]$  $[1 1 0 0 1 0]$ $[1 1 0 1 0 0]$ $[1 1 1 1 1 0]$   \\  \hline

& coext 14 &   $\bf [1 0 1 0 1 1]$  $\bf [1 1 1 0 1 0]$  \\  \hline

\hline
\end{tabular}
\end{center}
\normalsize
 \begin{center} Table 3: Cosimple single-element coextensions of $A$ $B$ and $C$ \end{center} 

By the Strong Splitter Theorem, the only cosimple single-element coextensions of $D$, $E$, and $F$ we must consider are the ones with $[0 0 0 0 0 1 1]$ as the new row. Let us call them $D'$, $E'$, $F'$, and $G'$, respectively. In each case we can find an $E_4$ minor. In particular, 
$D'/1\backslash \{3, 11\} \cong E_4$,  
$E'/1\backslash \{7, 11\} \cong E_4$, 
$F'/1\backslash \{3, 11\} \cong E_4$, and
$G'/1\backslash \{7, 11\} \cong E_4$.  $\qed$
\bigskip

Returning to the proof of Theorem 1.1, it remains to show that if $M$ has a $D_2$-minor and no $E_4$-minor, then we do not get any new matroids other than those already found in Lemma 2.2. 

Suppose $M$ is a cosimple single-element coextension of $D_2$. From Table 4 we see that $M$ is isomorphic to $A$, $B$, $C$, or $Z$. A matrix representation for $Z$ is shown below:

\small
\[ 
Z=\left[ 
\begin{array}{c|cccccc}
&    0&1&1&1&1&1 \\
&    1&0&1&1&1&0 \\
I_5& 1&1&0&1&0&0 \\
&    1&1&1&1&0&1 \\
&    0&0&0&1&1&1
\end{array} 
\right] 
\]
\normalsize

\noindent Since $Z$ is formed by adding only one row to $D_2$ (namely $[0 0 0 1 1 1]$) any coextension of $Z$ will also be a coextension of $A$, $B$, and $C$.

Suppose $M$ is a single-element extension of $D_2$. From  Table 1a we see that that $D_2$ has two single-element extensions $X_1$ and $X_3$ shown below:

 \small
\[ 
X_1=\left[ 
\begin{array}{c|ccccccc}
&   0&1&1&1&1&1&1\\
I_4&1&0&1&1&1&0&0 \\
&   1&1&0&1&0&0&1 \\
&   1&1&1&1&0&1&0
\end{array} 
\right] 
X_3=\left[ 
\begin{array}{c|ccccccc}
&   0&1&1&1&1&1&0\\
I_4&1&0&1&1&1&0&0 \\
&   1&1&0&1&0&0&1 \\
&   1&1&1&1&0&1&1
\end{array} 
\right] 
\] 
\normalsize

\noindent By the Strong Splitter Theorem the only coextensions of $X_1$ and  $X_2$ we must check are the ones formed with $[0 0 0 0 0 0 1 1]$ as the new row. Both these matroids have an $E_4$-minor. 

Lastly, suppose $M$ is a simple single-element extension of $Z$. It is straightforward to compute the three non-isomorphic simple single-element extensions which are $D$, $F$ and $Y$ (obtained by adding one of columns $[0 0 1 1 1]$, $[0 1 0 1 1]$, $[0 1 1 0 1]$, $[1 0 1 0 1]$, or $[1 1 1 0 0]$). The result follows again by Lemma 2.2 and the fact that when we add the above six columns to $Z$ we get the sixteen element matroid shown below which is isomorphic to $R_{17}\backslash \{17\}$.
\[
R_{16}=\left[ 
\begin{array}{c|cccccccccccc}
&    0&1&1&1&1&1&0&0&0&1&1 \\
&    1&0&1&1&1&0&0&1&1&0&1 \\
I_5& 1&1&0&1&0&0&1&0&1&1&1 \\
&    1&1&1&1&0&1&1&1&0&0&0 \\
&    0&0&0&1&1&1&1&1&1&1&0
\end{array}
\right] 
\]

\noindent Thus $Z$ does not contribute any new matroids to $EX[E_4]$ other than those found in Lemma 2.2. This completes the proof of Theorem 1.1. $\qed$

\small
 \begin{center}
\begin{tabular}{|c|p{25em}|c|c|}
\hline
Matroid &\bf{Coextension Rows} & {\bf Name} & Relevant minors  \\      \hline \hline
 
$D_2$ & $[0 0 0 0 1 1]$ $[0 0 0 1 0 1]$ $[0 0 0 1 1 0]$ $[0 0 1 1 1 1]$  $[1 0 0 1 1 1]$ $[1 0 1 0 0 0]$   &       $A_{26}$ ${\bf A}$          & $E_5$, $E_6^*$, $E_7$   \\  \hline

& $[0 0 0 1 1 1]$     &     $A_{31}$ ${\bf Z}$            &  $E_7$, $R_{10}$   \\  \hline

& $[0 0 1 0 0 1]$ $[0 1 0 1 0 0]$ $[0 1 1 1 0 1]$    & $A_{23}$  & $E_4$, $E_5$  \\  \hline

& $[0 0 1 0 1 0]$ $[0 0 1 1 0 0]$ $[0 1 0 0 0 1]$ $[0 1 0 0 1 0]$ $[0 1 1 0 1 1]$ $[0 1 1 1 1 0]$ &     $A_{20}$            & $E_4$, $E_6$    \\  \hline

& $[0 0 1 0 1 1]$ $[0 0 1 1 0 1]$  $[0 1 0 1 0 1]$ $[0 1 0 1 1 0]$ $[0 1 1 0 0 1]$ $[0 1 1 1 0 0]$ &  $A_{21}$   & $E_4$, $E_5$    \\  \hline

& $[0 0 1 1 1 0]$ $[0 1 0 0 1 1]$  $[0 1 1 0 1 0]$  &    $A_{24}$  & $E_4$    \\  \hline

& $[1 0 0 0 0 1]$ $[1 0 1 0 0 0]$  $[1 0 1 0 1 1]$ $[1 0 1 1 0 1]$ $[1 1 0 1 1 0]$ $[1 1 1 0 0 1]$  &  $A_{15}$   & $E_2$, $E_5$    \\  \hline

& $[1 0 0 0 1 0]$ $[1 0 0 1 0 0]$  $[1 1 0 0 0 0]$ $[1 1 0 1 0 1]$ $[11 1 1 0 0]$ $[1 1 1 1 1 1]$  &  $A_{6}$   & $E_1$, $E_4$    \\  \hline

& $[1 0 0 0 1 1]$ $[1 0 0 1 0 1]$  $[1 1 0 0 1 0]$ $[1 1 0 1 1 1]$ $[1 1 1 0 0 0]$ $[1 1 1 0 1 1]$ &    $A_{16}$  & $E_2$, $E_3$, $E_4$, $E_6^*$    \\  \hline

& $[1 0 0 1 1 0]$ $[1 0 1 0 1 0]$  $[1 0 1 1 0 0]$ $[1 0 1 1 1 1]$ $[1 1 0 0 0 1]$ $[1 1 1 1 1 0]$ &    $A_7$  & $E_4$, $E_5$   \\  \hline

& $[1 0 0 1 1 1]$ $[1 1 0 0 1 1]$  $[1 1 1 0 1 0]$   &    $A_{18}$ ${\bf C}$ & $E_3$, $E_5$, $E_6^*$, $E_7$    \\  \hline

& $[1 0 1 0 0 1]$  &    $A_{27}$ $\bf B$  & $E_5$    \\  \hline \hline

\end{tabular}
 \end{center}
\normalsize
 \begin{center} Table 4: Cosimple single-element coextensions of $D_2$ \end{center}

\bigskip
\noindent {\bf Proof of Corollary 1.2.}  Suppose $M$ is a $3$-connected binary  non-regular matroid with a $P_9$-minor, but no $P_9^*$-minor.  Lemma 2.1 implies that $P_9$ is a 3-decomposer or $M$ has a minor isomorphic to $D_2$ or $D_2^*$ since $E_4$ and $E_5$ have a $P_9^*$-minor. Observe that $D_2$ has four non-isomorphic cosimple single-element coextensions $A$, $B$, $C$, and $Z$. Of these $A$, $B$, and $C$ have an $E_5$-minor, and therefore a $P_9^*$-minor. Thus, if $M$ is a cosimple single-element coextension of $D_2$, then $M\cong Z$. Suppose $M$ is a simple single-element extension of $Z$, then $M\cong R_{16}$ as explained in the proof of Theorem 1.1. $\qed$
\bigskip

By duality we also have a decomposition result for $EX[P_9]$. Note that for $EX[P_9]$ among the rank 4 matroids only $F_7^*$, $S_8$, $AG(3,2)$ or $Z_4$ have no $P_9$-minor. All the others have a $P_9$-minor. Moreover, $S_8\cong Z_4\backslash b_r$ and $AG(3, 2)\cong Z_4\backslash c_r$, so we do not have to explicitly list them in the statement of the corollary.
\bigskip

Once the 3-decomposers and internally 4-connected extremal matroids are known for an excluded minor class, it is easy to determine all the internally 4-connected members of the class. Observe that $F_7$ and $F_7^*$ are internally 4-connected by default. All the restrictions of $PG(3,2)$ of size 10 and higher are internally 4-connected with the exception of $D_1$, $D_3$ and $X_2$ (8 matroids in total). The restrictions of $R_{16}$ in $EX[P_9^*]$ that are internally 4-connected are $R_{16}\backslash \{16\}$, $R_{16}\backslash \{15, 16\}$, $R_{16}\backslash \{14, 15, 16\}$, $R_{16}\backslash \{13, 14, 15, 16\}$, $R_{16}\backslash \{12, 13, 14, 15, 16\}$, and $R_{10}$ (six matroids in total). Thus there are 16 binary internally 4-connected matroids in $EX[P_9^*]$ and 8 in $EX[P_9]$.
\bigskip

We end this paper with a brief description of the usefulness of the Strong Splitter Theorem. The reader will note that it was used twice in the proof of Theorem 1.1. As a consequence we had to do calculations up to only rank 6 and 13 elements and that too for very  few matroids. Moreover, if needed, the reader can easily verify all the calculations. This is because once a required minor is found, checking that the specified deletions and contractions do indeed give the minor is very easy. Similarly for isomorphism. Once the isomorphism is found, it is easy to verify. The proof in this paper is not a ``computer-proof."

In [\ref{MayhewRoyleWhittle2011}] and [\ref{MayhewRoyle2012}] the authors used Chun, Mayhew, and Oxley's chain theorem for internally 4-connected binary matroids [\ref{ChunMayhewOxley2011}]. According to this chain theorem if a counterexample exists, then it has to be at most three more elements than one of the known internally 4-connected matroids in the class. So, for example, in [\ref{MayhewRoyle2012}] the authors extend and coextend known internally 4-connected matroids up to 3 more elements to see if any new internally 4-connected matroids appear. They had to do an exhaustive search up to rank 8 and 20 elements. Their proof is longer and had a component that relied soley on a computer check.

Had we not used the Strong Splitter Theorem and instead used Chun, Mayhew, and Oxley's result from [\ref{ChunMayhewOxley2011}],  we also would have had do an exhaustive search up to rank 8 and 20 elements. What's worse we wouldn't have found the splitter $M_{12}$ because it is not internally 4-connected, yet it has no 3-separation induced by the one in $P_9$ or $P_9^*$. That makes it an important member of the class. Due to the efficient stair-stepping approach of the Strong Splitter Theorem,  we had to coextend only four rank-5 matroids by 1 element (so the rank is at most 6) and extend by 2 elements up to 13 elements. Moreover, only four rank-6, 13-element matroids had to be checked for the required minor. 

Using the Strong Splitter Theorem gives an enormous savings in terms of computation, explains why an excluded minor class like $EX[E_4]$ lends itself to a concise characterization, and most importantly gives a short proof of a decomposition result which is stronger than a result that identifies all the internally 4-connected matroids. A complete list of all the internally 4-connected matroids follows immediately from a decomposition theorem like Theorem 1.1.

\bigskip


\noindent {\bf References}

\bigskip
\begin{enumerate}

\item  \label{ChunMayhewOxley2011} Chun, C., Mayhew, D. and Oxley, J. (2011) A chain theorem for internally 4-connected binary
matroids. {\it J. Combin. Theory Ser. B} 101, no. 3, 141 - 189.

\item  \label{KinganLemos2014}	Kingan, S. R.  and Lemos, M. (2014). Strong Splitter Theorem, {\it Annals of Combinatorics} {\bf 18-1}, 111-116.

\item \label{KinganLemos2014} Kingan, S. R.  and Lemos, M. (to appear). Decomposition of binary matroids with no prism minor, {\it Graphs and Combinatorics}.

\item \label{MayhewRoyleWhittle2011} Mayhew, D., Royle, G. and Whittle, G (2011). {\it The internally $4$-connected binary matroids with no $M(K_{3,3})$-minor},  Memoirs of the American Mathematical Society, {\bf 981}, American Mathematical Society, Providence, Rhode Island.

\item \label{MayhewRoyle2012} D. Mayhew and G. Royle (2012). The internally $4$-connected binary matroids with no $M(K_5\backslash e)$-minor, {\it Siam Journal on Discrete Mathematics}, {\bf 26}, 755-767.

\item \label{Oxley1987} Oxley, J. G. (1987) The binary matroids with no 4-wheel minor, {\it Trans. Amer. Math. Soc.}  {\bf 154}, 63-75.

\item \label{Oxley2012} Oxley, J. G. (2012) {\it Matroid Theory}, Second Edition, Oxford University Press, New York. 

\item  \label{Seymour1980} Seymour, P. D. (1980) Decomposition of regular matroids, {\it J. Combin. Theory  Ser. B} {\bf 28}, (1980) 305-359.

\end{enumerate}


\end {document}